\documentclass[pdflatex,sn-mathphys-num]{sn-jnl}

 \usepackage{enumerate}
\usepackage{graphicx}%
\usepackage{multirow}%
\usepackage{amsmath,amssymb,amsfonts}%
\usepackage{amsthm}%
\usepackage{mathrsfs}%
\usepackage[title]{appendix}%
\usepackage{xcolor}%
\usepackage{textcomp}%
\usepackage{manyfoot}%
\usepackage{booktabs}%
\usepackage{algorithm}%
\usepackage{algorithmicx}%
\usepackage{algpseudocode}%
\usepackage{listings}%

\def\f{\longrightarrow}
\def\im{\Longrightarrow}

\def\N{\mathbb{N}}

\def\e{\varepsilon}

\def\<{\langle}
\def\>{\rangle}

\def\A{\mathcal{A}}

\def\e{\varepsilon}

\def\R{\mathbb{R}}

\def\inte{\textnormal{int}\,}

\def\bdry{\textnormal{bdry}\,}

\def\nvarphi{\bp\nvarphi}
\def\sp{\hspace{0.015cm}}
\def\bp{\hspace{-0.08cm}}

\def\far{\textnormal{far}}
\makeatletter
\newcommand*\rel@kern[1]{\kern#1\dimexpr\macc@kerna}
\newcommand*\widebar[1]{%
  \begingroup
  \def\mathaccent##1##2{%
    \rel@kern{0.8}%
    \overline{\rel@kern{-0.8}\macc@nucleus\rel@kern{0.2}}%
    \rel@kern{-0.2}%
  }%
  \macc@depth\@ne
  \let\math@bgroup\@empty \let\math@egroup\macc@set@skewchar
  \mathsurround\z@ \frozen@everymath{\mathgroup\macc@group\relax}%
  \macc@set@skewchar\relax
  \let\mathaccentV\macc@nested@a
  \macc@nested@a\relax111{#1}%
  \endgroup
}
\makeatother

\theoremstyle{thmstyleone}%
\newtheorem{theorem}{Theorem}
\newtheorem{proposition}[theorem]{Proposition}%

\theoremstyle{thmstyletwo}%
\newtheorem{remark}{Remark}%

\theoremstyle{thmstylethree}%

\begin{document}

\title[Spherically Supported Property and the Exterior Sphere Condition with Infinite Radius]{Spherically Supported Property and the Exterior Sphere Condition with Infinite Radius}


\author*[1]{\fnm{Chadi} \sur{Nour}}\email{cnour@lau.edu.lb}

\affil*[1]{\orgdiv{Department of Computer Science and Mathematics}, \orgname{Lebanese American University}, \orgaddress{\city{Byblos Campus}, \postcode{P.O. Box 36}, \state{Byblos}, \country{Lebanon}}}


\abstract{We present several new characterizations of the spherically supported geometric property introduced in \cite{JCA2026}, emphasizing its connection with the exterior sphere condition with infinite radius. Moreover, we strengthen and provide a more direct and simpler proof of the main result established in \cite{JCA2026}.}

\keywords{Spherically supported sets, Strong convexity, Exterior sphere condition, Prox-regularity, Proximal analysis}



\maketitle

\section{Introduction}\label{sec1} Let $S\subset\R^n$ be a nonempty and closed set and let $s\in\bdry S$, the boundary of $S$. For $r>0$, a nonzero vector $\zeta\in N_S^P(s)$, the {\it proximal normal cone} to $S$ at $s$, is said to be {\it realized by an $r$-sphere}, if \begin{equation*} \label{realized}\left\<\frac{\zeta}{\|\zeta\|},x-s\right\>\leq \frac{1}{2r}\|x-s\|^2,\;\forall x\in S\;\;\;\; \left[\hbox{or equivalently}\; B\left(s+r\frac{\zeta}{\|\zeta\|};r\right)\cap S=\emptyset\right],\end{equation*}
where $B(y;\rho)$ denotes the open ball of radius $\rho$ centered at $y$. The same vector is said to be {\it far realized by an $r$-sphere} if \begin{equation*} \label{farrealized}\left\<\frac{\zeta}{\|\zeta\|},x-s\right\>\leq -\frac{1}{2r}\|x-s\|^2,\;\forall x\in S\;\;\;\; \left[\hbox{or equivalently}\; S\subset\widebar{B}\left(s-r\frac{\zeta}{\|\zeta\|};r\right)\right],\end{equation*}
where $\widebar{B}(y;\rho)$ denotes the closed ball of radius $\rho$ centered at $y$. For a fixed $r>0$, we recall the following four geometric properties: 
\begin{itemize}
\item We say that $S$ is {\it $r$-prox-regular} if for all $s\in\bdry S$ and for every nonzero vector $\zeta\in N_S^P(s)$, the normal vector $\zeta$ is realized by an $r$-sphere. Note that prox-regulairy implies the nontriviality (that is, $\not=\{0\}$) of the proximal normal cone $N_S^P(s)$ for any boundary $s$.  For further background on prox-regularity and related notions such as {\it positive reach}, {\it proximal smoothness}, {\it $p$-convexity}, and {\it $\varphi_0$-convexity}, see \cite{canino,csw,cm,fed,prt,shapiro,thibault}.
\item We say that $S$ satisfies the {\it exterior $r$-sphere condition} if for all $s\in\bdry S$, there exists a nonzero vector $\zeta\in N_S^P(s)$ such that the normal vector $\zeta$ is realized by an $r$-sphere. The exterior $r$-sphere condition, when imposed on the closure of the complement of $S$, corresponds to the classical {\it interior $r$-sphere condition} in control theory. It plays a crucial role in establishing regularity properties of the {\it minimal time function}; see Cannarsa and Frankowska~\cite{cf} and Cannarsa and Sinestrari~\cite{cs1,cs2}. Note that in \cite{JCA2024} (see also \cite{Pos2025}), Nour and Takche introduced an extended version of the exterior sphere condition, and proved that its complement is nothing but the union of closed balls with common  radius.
\item We say that $S$ is {\it $r$-strongly convex} if for all $s\in\bdry S$ and for every nonzero vector $\zeta\in N_S^P(s)$, the normal vector $\zeta$ is far realized by an $r$-sphere. Note that $r$-strong convexity is equivalent of $S$ being the intersection of closed balls with radius $r$. For further details on strong convexity and its applications, the reader is referred to the introduction of \cite{Nacry2025}, which offers a thorough historical overview of this property; see also \cite{balashov,frankowska,Goncharov,weber}.
\item We say that $S$ is {\it $r$-spherically supported} if for all $s\in\bdry S$, there exists a nonzero vector $\zeta\in N_S^P(s)$ such that the normal vector $\zeta$ is far realized by an $r$-sphere. This property was apparently introduced for the first time in \cite{JCA2026}.
\end{itemize}

It follows directly from the above definitions that if $S$ is $r$-prox-regular, then it satisfies the exterior $r$-sphere condition. The converse is not valid in general (see \cite[Example~2.5]{JCA2009}), but it holds when $S$ is {\it epi-Lipschitz} with compact boundary; see \cite[Corollary~3.12]{JCA2009}. Recall that $S$ is said to be {\it epi-Lipschitz} (or {\it wedged}) if, for every boundary point $s$, the set $S$ can be represented locally around $s$, after applying an orthogonal transformation, as the epigraph of a Lipschitz continuous function. This geometric notion, introduced by Rockafellar in \cite{rock}, can also be characterized by the nonemptiness of the topological interior of the Clarke tangent cone, which is equivalent to the pointedness of the Clarke normal cone; see \cite{clsw,rock}. Similarly, if $S$ is $r$-strongly convex, then it is $r$-spherically supported. The converse does not hold in general, as can be easily seen by taking $S$ to be the closed unit circle in $\mathbb{R}^2$. In \cite{JCA2026}, Nour and Takche proved that if $S$ is $r$-spherically supported with {\it nonempty interior}, then $S$ is $r$-strongly convex; see \cite[Theorem~1.1]{JCA2026}. This generalizes a known result in the literature (see $(\hbox{ii}')$ following the proof of \cite[Theorem~1.2]{weber}; see also \cite[Proposition~2.4]{Azagra}), where the equivalence is established under the {\it convexity} assumption on $S$. Note that the proof provided for \cite[Theorem~1.1]{JCA2026} turned out to be rather involved, since the set $S$ is assumed to have only a nonempty interior; hence, advanced tools from nonsmooth analysis are required.

The purpose of this paper is to investigate the $r$-spherically supported property and its relationship with the exterior $\infty$-sphere condition (that is, the exterior $r$-sphere condition for all $r>0$). It is natural to link these two notions, since the exterior $\infty$-sphere condition is always satisfied by any $r$-spherically supported set. Indeed, one can easily see that if a nonzero vector $\zeta \in N_S^P(s)$ is far realized by an $r$-sphere, then it is also realized by a $\rho$-sphere for every $\rho>0$. After providing new analytical characterizations of the $r$-spherically supported property, we prove that sets satisfying this property are either $r$-strongly convex or contained in the boundary of an $r$-strongly convex set. Likewise, we show that a set satisfying the exterior $\infty$-sphere condition is either convex or contained in the boundary of a convex set. As a consequence, we obtain a slightly strengthened version of \cite[Theorem~1.1]{JCA2026}, together with a more direct and simpler proof.

In the next section, we present the notations and basic definitions. Section~\ref{MR} is devoted to the statement and proof of our main results.

\section{Notations and basic definitions} We denote by $\|\cdot\|$, $\langle \cdot , \cdot \rangle$, $B$, and $\overline{B}$ the Euclidean norm, the usual inner product, the open unit ball, and the closed unit ball, respectively. For $r>0$ and $x \in \mathbb{R}^n$, we set 
\[
B(x;r) := x + rB \quad \text{and} \quad \overline{B}(x;r) := x + r\overline{B}.
\]
For a set $S \subset \mathbb{R}^n$, we denote by $S^c$, $\operatorname{int} S$, $\operatorname{bdry} S$, and $\operatorname{cl} S$ the complement (with respect to $\mathbb{R}^n$), the interior, the boundary, and the closure of $S$, respectively. The closed (resp.\ open) segment joining two points $x$ and $y$ in $\mathbb{R}^n$ is denoted by $[x,y]$ (resp.\ $]x,y[$).  For $\Omega \subset \mathbb{R}^n$ open and $f \colon \Omega \to \mathbb{R} \cup \{-\infty, +\infty\}$ an extended real-valued function, we denote by $\operatorname{epi} f$ the epigraph of $f$. The distance from a point $x$ to a nonempty and closed set $S \subset \mathbb{R}^n$ is denoted by 
\[
d_S(x) := \inf_{s \in S} \|x - s\|.
\]
We denote by $\operatorname{proj}_S(x)$ the set of points in $S$ closest to $x$, that is,
\[
\operatorname{proj}_S(x) := \{ s \in S : \|x - s\| = d_S(x) \}.
\]
The farthest distance from a point $x$ to a nonempty and closed set $S \subset \mathbb{R}^n$ is denoted by 
\[
d^f_S(x) := \sup_{s \in S} \|x - s\|.
\]
We denote by $\operatorname{far}_S(x)$ the set of farthest points in $S$ from $x$, that is,
\[
\operatorname{far}_S(x) := \{ s \in S : \|x - s\| = d^f_S(x) \}.
\]

Now we recall several notions from nonsmooth analysis that will be used throughout the paper. Comprehensive accounts of these concepts can be found in the monographs \cite{clsw,mord,penot,rockwet,thibault}.  Let $S \subset \mathbb{R}^n$ be a nonempty and closed set and $s \in S$. The {\it proximal normal cone} to $S$ at $s$, denoted by $N_S^P(s)$, is defined as
\[
N_S^P(s) := \bigl\{ \zeta \in \mathbb{R}^n : \exists\, \sigma \ge 0 \text{ such that } 
\langle \zeta, x - s \rangle \le \sigma \|x - s\|^2 \text{ for all } x \in S \bigr\}.
\]
A key feature of this cone, established in \cite[Proposition~1.1.5]{clsw}, is its {\it local nature}:  
whenever two closed sets coincide in a neighborhood of a point $s$, their proximal normal cones at $s$ also coincide. If $S$ is convex, the proximal normal cone $N_S^P(s)$ reduces to the {\it normal cone of convex analysis}, written $N_S(s)$ and given by
\[
N_S(s) := \bigl\{ \zeta \in \mathbb{R}^n : \langle \zeta, x - s \rangle \le 0,\; \forall x \in S \bigr\}.
\]
In this setting, $N_S(s)$ is nontrivial for every boundary point $s \in \operatorname{bdry} S$. The next statements summarize some geometric properties that hold for any nonempty and closed set $S \subset \mathbb{R}^n$. If $s \in \operatorname{proj}_S(x)$ for some $x \in \mathbb{R}^n$, then $s$ necessarily lies on the boundary of $S$, and the vector $x - s$ belongs to the proximal normal cone $N_S^P(s)$. Moreover, for every $y \in S$ one has
\[
\langle x - s, y - s \rangle \le \tfrac{1}{2}\|y - s\|^2,
\]
and the projection remains constant along the segment joining $s$ to $x$, that is,
\[
\operatorname{proj}_S(s + t(x - s)) = \{s\}, \quad \forall\, t \in [0,1[.
\]

\medskip
\noindent
Similarly, if $s \in \operatorname{far}_S(x)$, then $s \in \operatorname{bdry} S$ and the opposite vector $s - x$ belongs to $N_S^P(s)$. In this case,
\begin{equation*} \label{farnormal-new}
\langle s - x, y - s \rangle \le -\tfrac{1}{2}\|y - s\|^2, \quad \forall\, y \in S,
\end{equation*}
and the farthest point remains fixed along the ray starting from $x$ in the direction $x - s$, i.e.,
\[
\operatorname{far}_S(x + t(x - s)) = \{s\}, \quad \forall\, t > 0.
\]
Finally, for nonempty and closed set $S \subset \mathbb{R}^n$, $s\in \bdry S$ and $r>0$, we denote by:
\begin{itemize}
\item $[N^P_S(s)]_r:=\{\zeta\in N_S^P(s) : \hbox{is unit and realized by an}\;r\hbox{-sphere}\},$\vspace{0.1cm}
\item $[N^P_S(s)]^f_r:=\{\zeta\in N_S^P(s) : \hbox{is unit and far realized by an}\;r\hbox{-sphere}\}.$
\end{itemize}

\section{Main results}\label{MR} This section establishes our main results. We first provide analytical characterizations of the $r$-spherically supported property (Proposition~\ref{prop1}), then study the exterior $\infty$-sphere condition and its convexity implications (Theorem~\ref{th1}), and finally present the corresponding results for $r$-spherically supported sets (Theorems~\ref{th2} and~\ref{th3}).

We begin with the following proposition, which offers new analytical characterizations of the $r$-spherically supported property.
\begin{proposition} \label{prop1} Let $S \subset \mathbb{R}^n$ be nonempty and closed, and for $r>0$. The following assertions are equivalent:
\begin{enumerate}[$(i)$]
\item $S$ is $r$-spherically supported.
\item $S$ is bounded and for all $s\in\bdry S$ there exists $\zeta_s\in N_S^P(s)$ unit such that for all $R>0$,
$$\<\zeta_s-\zeta_x,x-s\>\leq -\frac{r+R}{2rR}\|x-s\|^2,\,\;\forall (x,\zeta_x)\in \bdry S\times [N^P_S(x)]^f_R.$$
\item $S$ is bounded and for all $s\in\bdry S$ there exists $\zeta_s\in N_S^P(s)$ unit such that for all $R>0$,
$$\|s-x\|\leq \frac{2rR}{r+R}\|\zeta_s-\zeta_x\|,\,\;\forall (x,\zeta_x)\in \bdry S\times [N^P_S(x)]^f_R.$$
\item $S$ is bounded and for all $s\in\bdry S$ there exists $\zeta_s\in N_S^P(s)$ unit such that for all $R>0$,
$$\<\zeta_s-\zeta_x,x-s\>\geq -\frac{2rR}{r+R}\|\zeta_s-\zeta_x\|^2,\,\;\forall (x,\zeta_x)\in \bdry S\times [N^P_S(x)]^f_R.$$
\end{enumerate}
\end{proposition}
\begin{proof} $(i)\im (ii)$: Let $s\in\bdry S$. Then there exists $\zeta_s\in N_S^P(s)$ unit such that \begin{equation} \label{nour1}\left\<\zeta_s,y-s\right\>\leq -\frac{1}{2r}\|y-s\|^2,\;\forall y\in S\;\;\;\; \left[\hbox{or equivalently}\; S\subset\widebar{B}\left(s-r\zeta_s;r\right)\right].\end{equation}
This clearly yields that $S$ is bounded. Now let $(x,\zeta_x)\in \bdry S\times [N^P_S(x)]^f_R$. Then we have \[\left\<\zeta_x,y-x\right\>\leq -\frac{1}{2R}\|y-x\|^2,\;\forall y\in S.\]
Taking $y:=x$ in this latter and $y:=x$ in \eqref{nour1}, we deduce that  \[\left\<\zeta_x,s-x\right\>\leq -\frac{1}{2R}\|s-x\|^2\;\,\hbox{and}\;\,\left\<\zeta_s,x-s\right\>\leq -\frac{1}{2r}\|x-s\|^2.\]
Hence \[\<\zeta_x-\zeta_s,s-x\>\leq- \frac{r+R}{2rR}\|s-x\|^2.\vspace{0.3cm}\]
$(ii)\im (iii)$:  Let $s\in\bdry S$. Then there exists $\zeta_s\in N_S^P(s)$ unit such that for all $R>0$,
$$\<\zeta_s-\zeta_x,x-s\>\leq -\frac{r+R}{2rR}\|x-s\|^2,\,\;\forall (x,\zeta_x)\in \bdry S\times [N^P_S(x)]^f_R.$$
This gives using Cauchy–Schwarz inequality that \[-\|s-x\|\|\zeta_s-\zeta_x\|\leq \<\zeta_x-\zeta_s,s-x\>\leq- \frac{r+R}{2rR}\|s-x\|^2,\,\;\forall (x,\zeta_x)\in \bdry S\times [N^P_S(x)]^f_R. \]
Hence, $$\|s-x\|\left(\frac{r+R}{2rR}\|s-x\| - \|\zeta_s-\zeta_x\|\right)\leq 0,\,\;\forall (x,\zeta_x)\in \bdry S\times [N^P_S(x)]^f_R.$$
Now, since the inequality in $(iii)$ is trivial when $s = x$, we may assume that $s \ne x$, which yields that
\[\|s - x\| \le \frac{2rR}{r + R}\,\|\zeta_s - \zeta_x\|.\vspace{0.3cm}\]
$(iii)\im (i)$: Let $s\in\bdry S$. Then there exists $\zeta_s\in N_S^P(s)$ unit such that for all $R>0$,
$$\|s-x\|\leq \frac{2rR}{r+R}\|\zeta_s-\zeta_x\|,\,\;\forall (x,\zeta_x)\in \bdry S\times [N^P_S(x)]^f_R.$$
We claim that $$S\subset\widebar{B}\left(s-r\zeta_s;r\right).$$
If not then, since $A$ is bounded, there exist $x\in\bdry A$ and $R>r$ such that $$A\subset \widebar{B}\left(x-r\zeta_s;R\right)\;\,\hbox{and}\;\,\|s-r\zeta_s-x\|=R.$$
Hence, $x\in \far_S(x-r\zeta_s)$ which yields that the unit vector $\zeta_x:=\frac{x-s+r\zeta_s}{R}\in N_S^Px)$ and is far realized by an $R$-sphere. Then $$\|r\zeta_s-R\zeta_x\|=\|s-x\|\leq \frac{2rR}{r+R}\|\zeta_s-\zeta_x\|.$$
Taking the square of the latter inequality with $c:=\<\zeta_s,\zeta_s\>\in [-1,1]$, we obtain that $$r^2+R^2-2rRc\leq \frac{8r^2R^2}{(r+R)^2} (1-c). $$
This yields that $$c\geq \frac{t^4+2t^3-6t^2+2t+1}{2t(t-1)^2},$$
where $t:=\frac{R}{r}>1$. Hence, $$c-1\geq \frac{t^4+2t^3-6t^2+2t+1}{2t(t-1)^2}-1=\frac{(t^2-1)^2}{2t(t-1)^2}=\frac{(t+1)^2}{2t}>0,$$
which contradicts $c\in [-1,1]$. \end{proof}

\begin{remark} \label{rem1} We cannot remove the boundedness assumption on $S$ in $(ii)$–$(iii)$ of Proposition~\ref{prop1}. Indeed, when $S$ is unbounded, we have $[N^P_S(x)]^f_R=\emptyset$ for all $x\in\bdry S$; hence, the inequalities in $(ii)$–$(iii)$ are automatically satisfied for every set $S$ with $N^P_S(s)\neq\{0\}$ at any boundary point $s$. It is also worth noting that the implication $(ii)\Rightarrow(iii)$ follows directly from the Cauchy–Schwarz inequality and therefore does not require the boundedness of $S$; this assumption is only needed in the geometric steps $(i)\Rightarrow(ii)$ and $(iii)\Rightarrow(i)$. \end{remark}

As mentioned in the introduction, if $S$ is $r$-spherically supported, then it satisfies the exterior $\infty$-sphere condition. Studying this latter condition thus provides deeper insight into the $r$-spherically supported property. On the other hand, it is known that if $S$ is $\infty$-prox-regular (that is, $r$-prox-regular for every $r>0$), then $S$ is convex. Since prox-regularity and the exterior sphere condition are closely related (see \cite{JCA2009}), a natural question arises: what type of convexity can be derived for $S$ if it satisfies the exterior $\infty$-sphere condition? The following theorem answers this question. In fact, it shows that if $S$ satisfies the exterior $\infty$-sphere condition, then it is either convex or contained in the boundary of a closed convex set.

\begin{theorem}\label{th1}
Let $S \subset \mathbb{R}^n$ be a nonempty, closed set satisfying the exterior $\infty$-sphere condition. Then $S$ is convex, or there exists a nonempty, closed, and convex set $A \subset \mathbb{R}^n$ such that $S \subset \operatorname{bdry} A$. In fact, $S$ is convex if $\inte S \not=\emptyset$, and a subset of the boundary of nonempty, closed, and convex set if $\inte S=\emptyset$.
\end{theorem}

\begin{proof} Let $S \subset \mathbb{R}^n$ be a nonempty, closed set satisfying the exterior $\infty$-sphere condition. For $s\in\bdry S$, we claim the existence of $\zeta_s\in N_S^P(s)$ unit such that  $$\<\zeta_s,x-s\>\leq 0,\;\,\forall x\in S.$$
Indeed, as $S$ satisfies the exterior $\infty$-sphere condition, we have, for each $n\in\N$, the existence of $\zeta_n\in N_S^P(s)$ unit  such that $$\<\zeta_n,x-s\>\leq \frac{1}{2n}\|x-s\|^2,\;\,\forall x\in S.$$
Taking $n\f+\infty$ in this latter and using that $(\zeta_n)_n$ has a subsequence, we do not relabel, that converges to a unit vector $\zeta_s$, we conclude that $$\<\zeta_s,x-s\>\leq 0,\;\,\forall x\in S,$$ which also yields that $\zeta_s\in N_S^P(s)$.\vspace{0.1cm}\\
{\it Case 1\sp}: $\inte S\not=\emptyset$.\vspace{0.1cm}\\
Then there exist $s_0\in S$ and $\rho>0$ such that $B(s_0;\rho_0)\subset S$. We claim that $S$ is epi-Lipschitz. If not, then by \cite[Exercise 3.6.5\sp(a)]{clsw} there exists $s \in \bdry S$ such that for all $u \in \R^n$
and $\varepsilon > 0$, one can find $s' \in S \cap B(s; \varepsilon)$, 
$t \in [0,\varepsilon[$, and $v \in B(u;\varepsilon)$ with 
$s' + t v \notin S$. Then there exists $s_n \in S \cap B(s; \tfrac{1}{n})$,
$v_n \in B(s-s_0; \tfrac{1}{n})$, and $t_n \in [0, \tfrac{1}{n}[$
such that $s_n + t_n v_n \notin S$.
\[\| s_n + v_n - s_0 \| 
\le \| s_n - s \| + \| v_n - (s_0 - s) \| 
\le \tfrac{2}{n}.
\]
For $n$ large, we have 
$s_n + v_n \in B(s_0; \rho) \subset \operatorname{int} S$. Then for $n$ large, there exists $\tau_n\in ]t_n,1[$ such that $s_n+\tau_nv_n\in \bdry S$. Then for $n$ large, there exists $\zeta_n \in N_S^P(s_n + \tau_n v_n)$ unit such that 
\[
\langle \zeta_n, x - s_n - t_n v_n \rangle \le 0,\;\, \forall x \in S.
\]
Then 
\[
\langle \zeta_n, s_n - s_n - t_n v_n \rangle \le 0\;\;\hbox{and}\;\;
\langle \zeta_n, s_n + v_n - s_n - t_n v_n \rangle \le 0,
\]
so
\[\langle \zeta_n, -t_nv_n \rangle \le 0 \;\;\hbox{and}\;\;
\langle \zeta_n, (1 - t_n) v_n \rangle \le 0.
\]
Then 
\[
\langle \zeta_n, v_n \rangle \ge 0 \;\;\hbox{and}\;\;\langle \zeta_n, v_n \rangle \le 0,
\]
which implies that \begin{equation}\label{nour2} \langle \zeta_n, v_n \rangle = 0.\end{equation} Since $s_n + v_n \in B(s_0; \rho) \subset \operatorname{int} S$,
we have 
\[
\| s_n + v_n - s_0 \| < \rho,\,\;\hbox{and hence,}\;\,
\rho - \| s_n + v_n - s_0 \| > 0.
\]
Let 
\[
\varepsilon_n = \frac{\rho - \| s_n + v_n - s_0 \|}{2} > 0.
\]
Then
\[
\| s_n + v_n + \varepsilon_n \zeta_n - s_0 \|
\le 
\| s_n + v_n - s_0 \| + \varepsilon_n 
< \| s_n + v_n - s_0 \| + \rho - \| s_n + v_n - s_0 \|= \rho.
\]
Thus 
\[
s_n + v_n + \varepsilon_n \zeta_n \in B(s_0; \rho) \subset \operatorname{int} S.
\]
Then
\[
\langle \zeta_n, s_n + v_n + \varepsilon_n \zeta_n - s_n - t_n v_n \rangle \le 0,
\]
which gives, using \eqref{nour2}, that
\[
 \e_n=\langle\zeta_n, (1 - t_n) v_n + \varepsilon_n \zeta_n \rangle \le 0,
\]
a contradiction. Therefore, $S$ is epi-Lipschitz. Now, from \cite[Theorem 7 $\&$ Remark 12]{NA2010}, we deduce that $S$ is $\infty$-prox-regular, and hence, $S$ is convex.\vspace{0.1cm}\\
{\it Case 2\sp}: $\inte S=\emptyset$.\vspace{0.1cm}\\
We define for any $s\in \bdry S=S$, $$[N^P_S(s)]_\infty:=\bigcap_{r>0}[N^P_S(s)]_r=\{\zeta_s\in N^P_S(s) : \<\zeta_s,x-s\>\leq 0,\;\,\forall x\in S\}. $$
Note that $[N^P_S(s)]_\infty$ is nonempty, as was shown above at the beginning of the proof. We define $$\A:=\{(s,\zeta_s) : s\in S\;\hbox{and}\;\zeta_s\in [N^P_S(s)]_\infty\},\;\,\hbox{and}$$ $$A:=\bigcap_{(s,\zeta_s)\in\A}\{x\in \R^n: \<\zeta_s,x-s\>\leq 0\}.$$
Clearly we have:
\begin{itemize}
\item $A$ is closed and convex with $$ \inte A=\{x\in\R^n : f(x)<0\}\;\,\hbox{and}\,\;\bdry A=\{x\in\R^n : f(x)=0\}.$$
where $$f(x):=\sup_{(s,\zeta_s)\in\A} \<\zeta_s,x-s\>,\;\,\forall x\in\R^n.$$
\item $S\subset A$.
\end{itemize}
We claim that $S\subset \bdry A$. Indeed, for $s\in S$, we have $$0\geq f(s)\geq \<\zeta_{s},s-s\>=0.$$ 
Therefore, $s\in\bdry A$, and hence, $S\subset \bdry A$.
\end{proof}

In light of Theorem~\ref{th1}, it is natural to investigate whether a similar result holds for the $r$-spherically supported property. In this case, the support is given by spheres of a fixed radius $r$ that contain $S$, which is a stronger requirement than the exterior $\infty$-sphere condition, as explained in the introduction. Therefore, one may expect a stronger form of convexity to emerge. The next theorem shows that this is indeed the case: an $r$-spherically supported set is either $r$-strongly convex or is contained in the boundary of an $r$-strongly convex set.\vspace{0.1cm}\\

\begin{theorem}\label{th2}
Let $S \subset \mathbb{R}^n$ be nonempty, closed, and $r$-spherically supported for some $r>0$. Then $S$ is $r$-strongly convex, or there exists a nonempty, closed, and $r$-strongly convex set $A \subset \mathbb{R}^n$ such that $S \subset \operatorname{bdry} A$. In fact, $S$ is $r$-strongly convex if $\inte S \not=\emptyset$, and a subset of the boundary of nonempty, closed, and $r$-strongly convex set if $\inte S=\emptyset$.
\end{theorem}
\begin{proof} As $S$ is $r$-spherically supported, we have that $[N^P_S(s)]^f_r\not=\emptyset$ for all $s\in\bdry S$, $S$ is compact,  and $S$ satisfies the exterior $\infty$-sphere condition.\vspace{0.1cm}\\
{\it Case 1\sp}: $\inte S\not=\emptyset$.\vspace{0.1cm}\\
Then by Theorem \ref{th1}, $S$ is convex. Hence, by \cite[Theorem 2.2.6]{schneider} combined with the definition of $r$-spherically supported property, we have that \begin{eqnarray*}\bigcap_{\substack{s\in\bdry S\\\zeta_s\in N^P_S(s)]^f_r}} \widebar{B}(s-r\zeta_s;r)\supset S&=& \bigcap_{\substack{s\in\bdry S\\\zeta_s\in N^P_S(s)]^f_r}} \{x\in\R^n : \<\zeta_s,x-s\>\leq 0\}\\&\supset& \bigcap_{\substack{s\in\bdry S\\\zeta_s\in N^P_S(s)]^f_r}} \widebar{B}(s-r\zeta_s;r).\end{eqnarray*}
This yields that $$S=\bigcap_{\substack{s\in\bdry S\\\zeta_s\in N^P_S(s)]^f_r}} \widebar{B}(s-r\zeta_s;r),$$ and hence, $S$ is $r$-strongly convex.\vspace{0.1cm}\\
{\it Case 2\sp}: $\inte S=\emptyset$.\vspace{0.1cm}\\
We define $$\A^f:=\left\{s-r\zeta_s : s\in S\;\hbox{and}\;\zeta_s\in [N^P_S(s)]^f_r \right\},\;\hbox{and}$$ $$A:=\bigcap_{c\in\A^f} \widebar{B}(c;r).$$
Clearly we have:
\begin{itemize}
\item $\A^f$ is compact.
\item $A$ is compact and $r$-strongly convex with \begin{eqnarray*} \inte A=\{x\in\R^n : f(x)<r\}\;\,\hbox{and}\,\;\bdry A&=&\{x\in\R^n : f(x)=r\}\\ &=&A\cap\bigcup_{c\in\A^f}\{x\in\R^n : \|c-x\|=r\},\end{eqnarray*}
where $f(\cdot)$ is the Lipschitz continuous function defined by $$f(x):=d^f_{\A^f}(x)=\sup_{c\in\A^f}\|c-x\|,\;\,\forall x\in\R^n.$$
\item $S\subset A$.
\end{itemize}
We claim that $S\subset \bdry A$. Indeed, for $s\in S$ and for $\zeta_s\in [N^P_S(s)]^f_r$,  we have $s\in A$ and $$s\in\{x\in\R^n : \|s+r\zeta_s-x\|=r\}\subset \bigcup_{c\in\A^f}\{x\in\R^n : \|c-x\|=r\}.$$ Therefore, $s\in\bdry A$, and hence, $S\subset \bdry A$.
\end{proof} 

As a consequence of Theorem \ref{th2}, we obtain a strengthened version of \cite[Theorem 1.1]{JCA2026}. 
The first improvement concerns the statement of \cite[Theorem 1.1]{JCA2026}. Indeed, since $N_S^P(\cdot)\subset N_{\bdry S}^P(\cdot)$, one can easily show that if $S$ is $r$-spherically supported, then its boundary is also $r$-spherically supported. However, the converse is not necessarily true, as shown by the example $S:=B^c$. The second improvement concerns the proof of the the necessary condition, which is now simplified and more direct.

\begin{theorem} \label{th3} Let $S\subset\R^n$ be a nonempty and closed set not reduced to a singleton, and let $r>0$. Then $S$ is $r$-strongly convex if and only if $S$ is bounded and $\bdry S$ is $r$-spherically supported with $\inte S\not=\emptyset$.
\end{theorem} 

\begin{proof} The sufficient condition is immediate. Indeed, when $S$ is a nonempty, closed set not reduced to a singleton and $r$-strongly convex, then it is bounded and $r$-spherically supported with $\operatorname{int} S \neq \emptyset$. Moreover, as mentioned above, since $N_S^P(\cdot) \subset N_{\operatorname{bdry} S}^P(\cdot)$, we directly obtain that $\bdry S$ is $r$-spherically supported.

We proceed to prove the necessary condition. Since $\bdry S$ is $r$-spherically supported and $\inte(\bdry S)=\emptyset$, we get from Theorem \ref{th2}, the existence of an $r$-strongly convex set $A$ such that $\bdry S\subset \bdry A$. We claim that $S=A$. Indeed, let $s\in\inte S$. If $s\not\in A$, then as $A$ is convex and using Hahn–Banach separation theorem, there exists $\zeta$ unit such that $$\<\zeta,a\>\leq 0\;\,\hbox{and}\,\;A\subset \{x\in\R^n : \<\zeta,x\>\geq 0\}. $$
Since $S$ is bounded and $s\in\inte S$, there exists $t>0$ such that $$a-t\zeta\in\bdry A\subset \bdry S.$$
Then, $$0\leq \<\zeta,a-t\zeta\>= \<\zeta,a\>-t\leq -t<0,$$
which gives the desired contradiction. Therefore, $\inte S\subset A$ which yields that $$S=\bdry S\cup \inte S\subset \bdry A\cup A=A.$$
For the other inclusion, assume that $S\subsetneq A$. Then there exists $x\in A$ such that $x\not\in S$. As $\inte S\not=\emptyset$, we consider $y\in\inte S\subset \inte A$. Since $A$ is convex, we get that $[y,x[\subset \inte A$. On the other hand, having $x\not\in S$ and $y\in\inte A$ yield the existence of $z\in]x,y[$ such that $z\in \bdry S\subset \bdry A$. This gives a contradiction as $z\in ]x,y[\subset ]x,y]\subset \inte A$. Therefore, $S\subset A$, and hence, $A=S$. This shows that $S$ is $r$-strongly convex. \end{proof}



\end{document}